\documentclass[a4paper]{amsart}

\usepackage{amssymb,amsfonts}

\usepackage[active]{srcltx}
\usepackage{graphicx}
\usepackage{amsmath,bbm,color,enumerate}
\usepackage{graphicx}
\usepackage[UKenglish]{babel}
\usepackage{lmodern}
\usepackage{microtype}
\usepackage{comment}
\usepackage{hyperref}

\newcommand*{\N}{\mathbb{N}}
\newcommand*{\Z}{\mathbb{Z}}
\newcommand*{\I}{\mathbb{I}}

\newcommand{\sulut}[1]{\left( #1 \right)}
\newcommand{\joukko}[1]{\left\{ #1 \right\}}

\newcommand*{\E}{\mathbb{E}}
\newcommand*{\Prob}{\mathbb{P}}

\begin{document}

\title[Expected characteristics]{Expected characteristic in Tunnels \& Trolls character creation, with generalizations}

\author[T. Brander]{Tommi Brander}
\address{DTU Compute -- Department of Applied Mathematics and Computer Science, 
Asmussens Alle, Building 322, entrance east,
DK-2800 Kgs. Lyngby, Denmark}
\email{tommi.brander@ntnu.no}


\thanks{Part of the research was carried out at Department of Mathematics and Statistics, P.O.Box 35 (MaD) FI-40014 University of Jyv\"askyl\"a, Finland.
I would like to thank Eija Laukkarinen for her comments and the anonymous referee for their comments.}

\keywords{stochastic sum, roleplaying game, dice}


\begin{abstract}
In the roleplaying game Tunnels \& Trolls the characteristics of player characters are determined by rolling dice in the following manner: First, one rolls three dice and calculates their sum. If the three dice all give the same result, another three dice are rolled and added to the total. This is continued until the three dice no longer match.
We calculate the average result of the stochastic sum: $10 + 4/5$.

We also consider a generalized dice rolling scheme where we roll an arbitrary number of dice with arbitrary number of sides. 
This generalization is motivated by various exotic dice that are used in many roleplaying games.
We calculate the expectation, and how much it differs from the situation where we only roll the set of dice once, with no rerolling and adding. As the number of dice increases, or the number of sides the dice have increases, this difference approaches zero, unless there are two dice (with the number of sides increasing), in which case the difference approaches one.
\end{abstract}

\newtheorem{theorem}{Theorem}[section]
\newtheorem{corollary}[theorem]{Corollary}
\newtheorem{lemma}[theorem]{Lemma}
\newtheorem{proposition}[theorem]{Proposition}

\theoremstyle{definition}
\newtheorem{definition}[theorem]{Definition}
\newtheorem{problem}[theorem]{Problem}
\newtheorem{example}[theorem]{Example}
\newtheorem{remark}[theorem]{Remark}

\numberwithin{equation}{section}

\maketitle


\section{Introduction}

Tunnels \& Trolls~\cite{StAndre:1975} is one of the oldest published tabletop roleplaying games~\cite[footnote~814]{Peterson:2012}.
Roleplaying games have become the subject of scholarly studies\footnote{In for example International journal of roleplaying~\url{http://journalofroleplaying.org/} and Analog game studies~\url{http://analoggamestudies.org/}.}, as the influence and study of video games has increased.
Many tabletop roleplaying games use dice to create or represent uncertainty~\cite{Dormans:2006,Torner:2014}.
Some of the schemes for rolling dice are mathematically nontrivial; we investigate one such scheme in the present article.
In tabletop roleplaying games, most players take the role of fictional characters. 
On defining roleplaying games see Arjoranta~\cite{Arjoranta:2011}, and for description see Dormans~\cite{Dormans:2006}.

In Tunnels \& Trolls, the player characters are described numerically by characteristics, which define for example how strong, lucky and intelligent the character is.
In the 5.5th and further editions of Tunnels \& Trolls~\cite{StAndre:2005} the characteristics are determined by rolling dice with a process explained below.

The starting characteristics are random, so we can hope to determine their expected value.
This value is useful for a player of the game, since it allows quickly judging the worth of a new player character -- are its scores above or below average?
A designer or game master might also consider various alternative methods of determining characteristics, and information about their average is helpful in the process.
The expected values and distributions of various methods of rolling dice are often discussed by game designers and players of tabletop roleplaying games.

Each characteristic is independently determined with the following process.
\begin{enumerate}
\item Roll three dice (ordinary six-sided dice).
\item If the result is triples -- each of the dice shows the same value -- roll three additional dice.
\item Continue until the newly rolled dice are not a triple.
\item Sum the results of all the dice rolled thus far.
\end{enumerate}

\begin{example}
Suppose we roll $(3, 3, 1)$. This is not a triple, so the value of the characteristic is $3+3+1 = 7$.
\end{example}

\begin{example}
Suppose we roll $(5, 5, 5)$. This is a triple, so we roll again and get $(1,1,1)$. This is also a triple, so we roll yet again and get $(1, 2, 2)$. This is no longer a triple, so we sum everything rolled thus far: $5+5+5+1+1+1+1+2+2 = 23$.
\end{example}

The same rolling process can be used with an arbitrary number of dice, which may have an arbitrary number of sides.
Indeed, the saving roll system in Tunnels \& Trolls uses the same process with two six-sided dice; there, doubles add and roll over.
It turns out that the expected value of a characteristic is $54/5 = 10.8$, and the expected value of the dice roll in a saving roll is $42/5 = 8.4$.
These values can be found as corollary~\ref{cor:expect_tt}.

Several roleplaying games use various exotic dice (for example Pathfinder~\cite{Bulmahn:2009} and Dungeon Crawl Classics~\cite{Curtis:Goodman:Stroh:Zimmerman:2012}), which motivates us to ask how the number of sides the dice have influences the expectation.

It is interesting to ask how significant it is to roll all the dice again when they match and add this to the previous result.
One way of answering this is checking how much the expected value of the result changes when matches are added and rolled over, when compared to the situation where only one set of dice is rolled and matching dice do not have any special meaning.
It turns out that the difference in expectations vanishes in the limit of increasing number of dice or increasing number of sides in the dice, except when we roll precisely two dice, and let the number of sides they have increase to infinity.
For the precise result, see theorem~\ref{thm:expect_diff}.

\section{Formalization and calculation}

We assume we are rolling $n \in \N$ dice, each of which has $s \in \Z_+$ sides.
Let both $n$ and $s$ be fixed.
The cases $s = 1$ or $n \in \joukko{0,1}$ are trivial, so by default we suppose $n,s \geq 2$.
See remark~\ref{remark:trivial} for a precise formulation of the triviality.

Each $s$-sided die is represented by a random variable $Z^l_j$ where the index $j \in \Z_+$ indicates which set of dice rolls is in question and $l \in \joukko{1,2,\ldots,n}$ indicates which of the dice in a given set of dice rolls is in question; this is made precise in equation~\eqref{eq:x} and the text below it.
Each such random variable is a mapping $Z^l_j \colon \Omega \to \N$ such that for every $ m \in \joukko{1,\ldots,s}$ we have $\Prob \joukko{\omega \in \Omega; Z^l_j(\omega) = m} = 1/s$, and probability of every other event is zero.
We assume the die rolls, i.e.\ variables $Z^l_j$, are mutually independent.

Let $M_0 = \Omega$ and $M_j$ be the event of the $j$th set of rolls being a match; that is, for $j \ge 1$,
\begin{equation}
M_j = \joukko{\omega \in \Omega; Z_j^1 = Z_j^2 = \cdots = Z_j^n}.
\end{equation}

Define the random variables $X_j$ as follows: $X_0 = 0$ and for  $j \ge 1$
\begin{equation} \label{eq:x}
X_j =
\begin{cases}
X_{j-1} + \sum_{l=1}^s Z_j^l  \text{ when } \omega \in \bigcap_{k = 0}^{j-1} M_{k} \\
X_{j-1} \text{ otherwise.}
\end{cases}
\end{equation}
The interpretation of these random variables in terms of dice rolls is that $X_1$ is the sum of the first set of dice rolls, $X_2$ allows for the possibility of matches and the corresponding additional dice rolls, $X_3$ allows for two sets of matches, and so on.
Hence, we define
\begin{equation}
X^{n,s} = \lim_{j \to \infty} X_j,
\end{equation}
which is the final outcome of the entire dice rolling  scheme.
Note that the fixed constants $n$ and $s$ are implicit in the $X_j$ variables.
We sometimes write them explicitly as superindices.

\begin{remark}\label{remark:trivial}
We have $X^{0,s} = 0$ and $X^{1,s} = X^{n,1} = \infty$, where the free variables satisfy $n,s\geq 1$.
\end{remark}

\begin{remark}
Characteristics in Tunnels \& Trolls are rolled with $X^{3,6}$ and saving rolls with $X^{2,6}$.

In saving rolls there is an additional rule of automatic failure when the initial roll gives the result of three.
If one wanted to assign it a numerical value of $0$ (an arbitrary choice which would not necessarily lead to failure in game), then the relevant expectation would have to decreased by 
\begin{equation}
\Prob\sulut{X^{2,6}=3}\cdot 3 = \frac{2}{36} \cdot 3 = 1/6.
\end{equation}
This type of adjustment is quite easy to do and generalizing it to all of the situations covered in the paper is far from obvious, so we ignore it from now on.
We stress that the choice of zero here is both arbitrary and unsatisfactory; in theory, no finite numerical value leads to guaranteed failure, though large negative numbers such as $-1000$ would almost certainly guarantee failure in all practical game situations.
\end{remark}

We are interested in the expected value $\E \sulut{X^{n,s}}$, and the relation between it and $\E\sulut{X^{n,s}_1}$.
That is, we want the know how much the expected value increases when matches allow rolling an additional set of dice.

The following lemma follows from the linearity of expectation and the definition of the random variables:
\begin{lemma}[Expectation without rerolls]
For $n \in \N$ and $s \in \Z_+$ we have
\begin{equation}
\E \sulut{X^{n,s}_1} = n\sulut{s+1}/2.
\end{equation}
\end{lemma}

\begin{corollary}[Expectations of dice rolls in Tunnels \& Trolls without rerolls]
\begin{align}
\E\sulut{X^{3,6}_1} = 21/2  \text{ and }
\E\sulut{X^{2,6}_1} = 7.
\end{align}
\end{corollary}

By monotone convergence, $\E (X^{n,s}) = \lim_{j \to \infty} \E(X_j) = \E(X_N)$, where
\begin{equation}
N(\omega) = \sup\joukko{j \in \N; \omega \in \bigcap_{k = 0}^{j-1} M_{k}}
\end{equation}
is the first time the dice do not match.
Note that $N$ is finite almost surely.

We now calculate the probability mass function~$f$ of $N$.
The probability of rolling a match is $p = \sulut{\frac{1}{s}}^{n-1} = s^{1-n}$, so
\begin{equation}
\begin{split}
f(j) &=
\begin{cases}
p^{j-1}  \sulut{1-p} &\text{ for } j \geq 1 \\
0 &\text{ for } j = 0.
\end{cases} 
\end{split}
\end{equation}

On the other hand,
\begin{equation}
\E(X^{n,s}) = \E \sulut{\E\left[ X^{n,s} | N \right]} = \sum_{j=1}^\infty f(j) \E \left[ X_j | N = j \right].
\end{equation}
To calculate this explicitly we need to know the value of the conditional expectation $\E \left[ X_j | N = j \right]$.
For all $k < j$ we have $\omega \in M_k$, so $Z_k^1 = \ldots = Z_k^n$.
Further, the set of random variables $\joukko{N}\cup \joukko{Z_l^1; l \in \N}$ is independent.
Thus,
\begin{equation}
\begin{split}
\E \left[ X_j | N = j \right] &= \E\left[\sum_{l=1}^n Z_j^l + \sum_{k=1}^{j-1} \sum_{l=1}^n Z_k^l \; \Big| N = j\right] \\
&= \E\left[\sum_{l=1}^n Z_j^l\; \Big| N = j\right] + \sum_{k=1}^{j-1}\E\left[ n Z_k^1 \; \Big| N = j\right].
\end{split}
\end{equation}
We calculate the first and the second part of the expectation separately, starting from the second sum:
\begin{equation}
\begin{split}
\sum_{k=1}^{j-1}\E\left[ n Z_k^1 \; \Big| N = j\right] 
&= n\sum_{k=1}^{j-1}\E\sulut{ Z_k^1 } \\
&= n\sulut{j-1}\sulut{s+1}/2.
\end{split}
\end{equation}
For the first sum we have
\begin{equation}
\E\left[\sum_{l=1}^n Z_j^l\; \Big| N = j\right] = \E\left[\sum_{l=1}^n Z_j^l\; \Big| \omega \in \bigcap_{k=0}^{j-1} M_k \setminus M_j \right].
\end{equation}
By the independence of all the $Z$ variables, and in particular the independence of $Z_j$ from $M_k$ when $k < j$ and the independence of $M_k$ from $M_j$, we get
\begin{equation}
\begin{split}
\E\left[\sum_{l=1}^n Z_j^l\; \Big| \omega \in \bigcap_{k=0}^{j-1} M_k \setminus M_j \right] 
&= \E\left[\sum_{l=1}^n Z_j^l\; \Big| \omega \in \Omega \setminus M_j \right] \\
&= \E\sulut{\sum_{\omega \in \Omega \setminus M_j} \sum_{l=1}^n Z_j^l(\omega) } / \Prob\sulut{\Omega \setminus M_j} \\
&= \E\sulut{\sum_{\omega \in \Omega} \sum_{l=1}^n Z_j^l(\omega) - \sum_{\omega \in M_j} \sum_{l=1}^n Z_j^l(\omega)} / \sulut{1-p} \\
&= \sulut{\E\sulut{ \sum_{l=1}^n Z_j^l} - \E\sulut{ \sum_{l=1}^n Z_j^l \I_{M_j}}} / \sulut{1-p} \\
&= \sulut{n\sulut{s+1}/2 - n\E\sulut{ Z_j^1} \Prob\sulut{M_j}} / \sulut{1-p} \\
&=n\sulut{s+1}/2.
\end{split}
\end{equation}

Thus we have
\begin{equation}
\E \left[ X_j | N = j \right] =  j n(s+1)/2,
\end{equation}
whence
\begin{equation}
\E(X^{n,s}) = \sum_{j=1}^\infty p^{j-1}  \sulut{1-p} j n(s+1)/2.
\end{equation}
This is a series of the form $c\sum_{j=1}^\infty j a^{j-1}$,
which has the value
\begin{equation}
c\sum_{j=1}^\infty j a^{j-1} = c\sulut{1-a}^{-2}.
\end{equation}
So, we have
\begin{equation}
\begin{split}
\E\sulut{X^{n,s}}
&= \E\sulut{X_1^{n,s}} / (1-s^{1-n}).
\end{split}
\end{equation}
This proves the following theorem:
\begin{theorem}[Expectation of the random sum]
Suppose $n \geq 2$ and $s \geq 2$. Then 
\begin{equation}
\E\sulut{X^{n,s}} = (1-s^{1-n})^{-1} n(s+1)/2.
\end{equation}
\end{theorem}

\begin{corollary}[Expectations related to Tunnels \& Trolls] \label{cor:expect_tt}
\begin{align}
\E\sulut{X^{3,6}} = 54/5  \text{ and }
\E\sulut{X^{2,6}} = 42/5.
\end{align}
\end{corollary}

The identity in the theorem is consistent with the trivial identities $\E\sulut{X^{0,s}} = 0$, $\E\sulut{X^{1,s}} = \infty$ and $\E\sulut{X^{n,1}} = \infty$, and furthermore, as $n \to \infty$, we have
\begin{equation}
\begin{split}
\E\sulut{X^{n,s}}-\E\sulut{X^{n,s}_1}  &= \sulut{\frac{1}{1-s^{1-n}}-1}\frac{n(s+1)}{2} \\
&= \sulut{\frac{ns^{1-n}}{1-s^{1-n}}}\frac{s+1}{2} \to 0.
\end{split}
\end{equation}
On the other hand, as $s \to \infty$, we have
\begin{equation}
\begin{split}
\E\sulut{X^{n,s}}-\E\sulut{X^{n,s}_1}  &= \sulut{\frac{s^{2-n}+s^{1-n}}{1-s^{1-n}}}\frac{n}{2} \to
\begin{cases}
0 \text{ when } n \geq 3 \\
1 \text{ when } n = 2.
\end{cases}
\end{split}
\end{equation}
Thus, we get the following theorem:
\begin{theorem}[Limits of expectations]  \label{thm:expect_diff}
\begin{equation}
\begin{split}
\text{For } s\in \Z_+: \E\sulut{X^{0,s}}-\E\sulut{X^{0,s}_1}  &= 0. \\
\text{For } s\in \Z_+: \E\sulut{X^{1,s}}-\E\sulut{X^{1,s}_1}  &= \infty. \\
\text{For } n\in \Z_+: \E\sulut{X^{n,1}}-\E\sulut{X^{n,1}_1}  &= \infty. \\
\text{For } s \ge 2:\lim_{n \to \infty} \sulut{\E\sulut{X^{n,s}}-\E\sulut{X^{n,s}_1}}  &= 0. \\
\lim_{s \to \infty} \sulut{\E\sulut{X^{2,s}}-\E\sulut{X^{2,s}_1}}  &= 1. \\
\text{For } n \ge 3 : \lim_{s \to \infty} \sulut{\E\sulut{X^{n,s}}-\E\sulut{X^{n,s}_1}}  &= 0.
\end{split}
\end{equation}

\end{theorem}

\section{Conclusion}

The effect of rolling again and adding dice is quite small, unless one is rolling two fairly small dice.
The probability of very high results does grow from zero to a small but positive number, which might be relevant even in the absence of a large change in the expected value.


\bibliographystyle{plain}
\bibliography{math}

\begin{thebibliography}{1}

\bibitem{StAndre:1975}
Kenneth Eugene~St.\ Andre.
\newblock {\em Tunnels \& Trolls}.
\newblock April 1975.

\bibitem{StAndre:2005}
Kenneth Eugene~St.\ Andre.
\newblock {\em Tunnels \& Trolls}.
\newblock Flying Buffalo Incorporated, 5.5th edition, 2005.

\bibitem{Arjoranta:2011}
Jonne Arjoranta.
\newblock Defining role-playing games as language-games.
\newblock {\em International Journal of Role-Playing}, 2:3--17, 2011.

\bibitem{Bulmahn:2009}
Jason Bulmahn.
\newblock {\em Pathfinder Roleplaying Game Core Rulebook}.
\newblock Pathfinder RPG. Paizo Publishing, LLC, August 2009.
\newblock System reference document at
  \url{http://paizo.com/pathfinderRPG/prd/}.

\bibitem{Curtis:Goodman:Stroh:Zimmerman:2012}
Michael Curtis, Joseph Goodman, Harley Stroh, and Dieter Zimmerman.
\newblock {\em Dungeon Crawl Classics Role Playing Game}.
\newblock Goodman Games, 2012.

\bibitem{Dormans:2006}
Joris Dormans.
\newblock On the role of the die: A brief ludologic study of pen-and-paper
  roleplaying games and their rules.
\newblock {\em Game Studies}, 6(1), 2006.
\newblock Available \url{http://gamestudies.org/0601/articles/dormans}.

\bibitem{Peterson:2012}
Jon Peterson.
\newblock {\em Playing at the world}.
\newblock Unreason Press, 2012.

\bibitem{Torner:2014}
Evan Torner.
\newblock Uncertainty in analog roleplaying games, part 2.
\newblock {\em Analog game studies}, 1(2), September 2014.
\newblock Available
  \url{http://analoggamestudies.org/2014/09/uncertainty-in-analog-role-playing-games-part-2/}.

\end{thebibliography}

\end{document}